# B-splines and kernels
# of trigonometric interpolation splines

Volodymyr Denysiuk

Dr. of Ps.-M. Sciences, Professor, Kiev, Ukraine

kvomden@nau.edu.ua

Olena Hryshko

Senior Lecturer, National Aviation University, Kiev, Ukraine

olena.hryshko@gmail.com

**Abstracts.**

This article focuses on trigonometric Riemann B-splines and Riemann kernels of trigonometric interpolation splines of arbitrary order; it is shown that trigonometric interpolation splines are a convolution of trigonometric B-splines with corresponding kernels. Theoretical statements are followed by examples, and the results are applicable in many practical areas.

**Introduction.**

The approximation and representation of an arbitrary known or unknown function by a set of some special functions is a central topic of our analysis. The term "special functions" refers to classes of algebraic and trigonometric polynomials and their modifications; at the same time, trigonometric series are considered included in the trigonometric polynomials classes. Usually, such special functions are easy to compute and have interesting analytical properties [1].

One of the most successful modifications of algebraic polynomials is polynomial splines, which are sewn from segments of these polynomials according to a certain scheme. The theory of polynomial splines is quite recent yet well-developed (see, for instance, [2], [3], [4], [5]). One of the advantages of polynomial splines is their approximation properties [6]. The main disadvantage of polynomial splines is their piecewise structure, which greatly complicates their use in analytical transformations.

An important role in the theory of polynomial splines is played by normalized B-splines (hereinafter simply B-splines), which are a basis in the spline function space [3], [4]. The representation of polynomial splines through B-splines is used in many problems requiring analytical transformations, in particular, in solving differential and integral equations (see, for instance, [7], [8]).

There are also modifications of trigonometric series [9], [10] that depend on several parameters and have the same properties as polynomial splines [11], [12]; moreover, the class of such modified series is quite wide and includes the class of polynomial periodic splines. As a result, the class of such series' is called trigonometric spline.

Trigonometric splines and their generalizations were also discussed in [13], and [14].

Trigonometric B-splines and their generalizations can play a special role in the theory of trigonometric splines, as well as in the theory of polynomial splines. Such splines were discussed in [15]. However, the approach used in this paper does not allow us to explore some important properties of B-splines.

**Purpose of the work.** Development of a method for trigonometric B-splines construction; exploration of their properties.

**The main part of the paper.**

Let us set a continuous periodic function $f(t)$ on the segment $[0, 2\pi]$. A uniform grid $\Delta_N = \{t_i\}_{i=1}^N$, $t_i = \frac{2\pi}{N}(i-1)$, $N = 2n+1$, $n = 1, 2, \ldots$ is also set on this segment. Let us also set $\lambda = \pi/N$. We denote by $\{f(t_i)\}_{i=1}^N = \{f_i\}_{i=1}^N$ the set of output values of the function $f(t)$ at the nodes of the grid $\Delta_N$. Let us have a look at the trigonometric polynomial

$$T_n^*(t) = \frac{a_0^*}{2} + \sum_{k=1}^n a_k^* \cos kt + b_k^* \sin kt, \quad (1)$$

that interpolates the function $f(t)$ on the grid $\Delta_N$. Then the coefficients of this polynomial are determined by the formulas

$$a_k^* = \frac{2}{N} \sum_{j=1}^N f_j \cos kt_j, \qquad b_k^* = \frac{2}{N} \sum_{j=1}^N f_j \sin kt_j,$$
$$k = 0, 1, \ldots, n; \qquad\qquad k = 1, 2, \ldots, n. \quad (2)$$

Let us set $C_p^r$ to be a class of $2\pi$ periodic continuous functions having continuous derivatives up to and including $r-1$ order. Let $C_p^0$ be a class of $2\pi$ periodic continuous functions, and $C_p^{-1}$ be a class of piecewise continuous functions with a finite number of jump discontinuities.

In [9], [10], it was shown that a trigonometric spline $St(f, \sigma, r, N, t) \in C^{r-2}$, interpolating a function $f(t)$ at the nodes of the grid $\Delta_N$, can be represented as

$$St(f, \sigma, r, \Delta_N, t) = \frac{a_0}{2} + \sum_{k=1}^{\frac{N-1}{2}} \left[ H_k^{-1}(\sigma, r, k) \left( a_k C_k(\sigma, r, N, t) + b_k S_k(\sigma, r, N, t) \right) \right], \quad (3)$$

where $\sigma_k(r, N) = \left( \frac{\sin(k\pi/N)}{k} \right)^{1+r}$ is a Riemann multiplier, which is included in the generalized summation method (R,r) [16];

$H_k(\sigma, r, k)$ - interpolation multiplier,

$$H_k(\sigma, r, k) = \sigma(r, k, N) + \sum_{m=1}^{\infty} \left[ \sigma(r, k, mN - k) + \sigma(r, k, mN + k) \right]; \quad (4)$$

$C_k(\sigma, r, N, t)$, $S_k(\sigma, r, N, t)$ - spline functions defined as follows:

$$C_k(\sigma, r, N, t) = \sigma(r, k, N) \cos(kt) +$$
$$+ \sum_{m=1}^{\infty} \left[ \sigma(r, k, mN - k) \cos((mN - k)t) + \sigma(r, k, mN + k) \cos((mN + k)t) \right]; \quad (5)$$

$$S_k(\sigma, r, N, t) = \sigma(r, k, N) \sin(kt) +$$
$$+ \sum_{m=1}^{\infty} \left[ -\sigma(r, k, mN - k) \sin((mN - k)t) + \sigma(r, k, mN + k) \sin((mN + k)t) \right],$$
$$(k = 1, 2, \ldots, n).$$

There is an equality for all positive integers $k, m, n$:

$$\sigma_k(m + n, N) = \sigma_k(m, N) \sigma_k(n, N). \quad (6)$$

Indeed,

$$\sigma_k(m+n,N) = \left(\frac{\sin(k\pi/N)}{k}\right)^{1+(m+n)} = \left(\frac{\sin(k\pi/N)}{k}\right)^{1+m} \left(\frac{\sin(k\pi/N)}{k}\right)^{1+(n-1)} = \sigma_k(m,N)\sigma_k(n-1,N).$$

Expression (3) is a trigonometric Fourier series that converges uniformly to $r > 1$, with coefficients

$$A_{mN\pm k} = a_{mN\pm k} H_k^{-1}(r,N)\sigma_{mN\pm k}(r,N),$$
$$B_{mN\pm k} = b_{mN\pm k} H_k^{-1}(r,N)\sigma_{mN\pm k}(r,N),$$
$$(k=1,2,...,n;\ m=0,1,...).$$

Taking into account (6), we can represent these coefficients the following way:

$$A_{mn\pm k} = \left[a_{mN\pm k} H_k^{-1}(r,N)\sigma_{mN\pm k}(0,N)\right]\sigma_{mN\pm k}(r,N),$$
$$B_{mn\pm k} = \left[b_{mN\pm k} H_k^{-1}(r,N)\sigma_{mN\pm k}(0,N)\right]\sigma_{mN\pm k}(r,N) \tag{7}$$
$$(k=1,2,...,n;\ m=0,1,...).$$

The coefficients of the trigonometric series (5) are the products of two expressions. According to the general theory of trigonometric series [17] [18], the sum of such a series is the convolution of two functions. These two functions are represented by the following Fourier series:

$$BR(r,t) = \frac{1}{\pi}\left[\frac{1}{2} + \sum_{k=1}^{\frac{N-1}{2}} C_k(\sigma,r,N,t)\right] \tag{8}$$

and

$$KR0(f,\sigma,2j,N,t) = \frac{a_0}{2} + \sum_{k=1}^{\frac{N-1}{2}}\left[\frac{C0_k(\sigma,0,N,t)}{Hc(2j,k)}a_k + \frac{S0_k(\sigma,0,N,t)}{Hs(2j,k)}b_k\right], \tag{9}$$

$$KR1(f,\sigma,2j-1,N,t) = \frac{a_0}{2} + \sum_{k=1}^{\frac{N-1}{2}}\left[\frac{C1_k(\sigma,0,N,t)}{Hc(2j-1,k)}a_k + \frac{S1_k(\sigma,0,N,t)}{Hs(2j-1,k)}b_k\right], \tag{9a}$$

$$(j=1,2,\ ...\ ),$$

where

$$C0_k(\sigma,0,N,t) = \sigma(0,k)\cos(kt) +$$
$$+ \sum_{m=1}^{\infty}(-1)^m\left[\sigma(0,2mN+k)\cos((2mN+k)t) + \sigma(0,2mN-k)\cos((2mN-k)t)\right];$$

$$S0_k(\sigma,r,N,t) = \sigma(0,k)\sin(kt) +$$
$$+ \sum_{m=1}^{\infty}(-1)^m\left[\sigma(0,2mN+k)\sin((2mN+k)t) - \sigma(0,2mN-k)\sin((2mN-k)t)\right].$$

and

$$C1_k(\sigma,0,N,t) = \sigma(0,k)\cos(kt) +$$
$$+ \sum_{m=1}^{\infty}\left[\sigma(0,2mN+k)\cos((2mN+k)t) + \sigma(0,2mN-k)\cos((2mN-k)t)\right];$$

$$S1_k(\sigma,r,N,t) = \sigma(0,k)\sin(kt) +$$
$$+ \sum_{m=1}^{\infty}\left[\sigma(0,2mN+k)\sin((2mN+k)t) - \sigma(0,2mN-k)\sin((2mN-k)t)\right]$$

The function $BR(r,t)$ is called a trigonometric Riemann B-spline of order $r$ ($r = 0,1,...$), and the function $KR0(f,\sigma,2j,N,t)$ and $KR1(f,\sigma,2j+1,N,t)$ is called the Riemann kernel of the interpolation trigonometric spline $St0(f,\sigma,2j,N,t)$ and $St0(f,\sigma,2j-1,N,t)$. This interpolation trigonometric spline $St(f,\sigma,r,N,t)$ can be represented as a convolution:

$$St0(f,\sigma,2j,N,t) = \int_0^{2\pi} KR0(f,\sigma,2j,N,v) BR(2j-1,t-v) dt, \qquad (10)$$

$$St1(f,\sigma,2j-1,N,t) = \int_0^{2\pi} KR1(f,\sigma,2j-1,N,v) BR(2(j-1),t-v) dt, \qquad (10a)$$

$$(j = 1, 2, \ldots).$$

For illustration, we provide the plots of the first few values of the trigonometric B-splines and Riemann kernels for the interpolation spline $St(f,\sigma,r,N,t)$ interpolating the trigonometric polynomial (1) at the grid $\Delta_N$ nodes, assuming $N=9$, and $\{f_i\}_{i=1}^{9} = \{2,1,3,2,4,1,3,1,3\}$. Note the abbreviated notations used in the plots: the trigonometric interpolation spline $St(f,\sigma,r,N,t)$ is denoted as $St(r,t)$, and the kernel of this spline $KR0(f,\sigma,2j,N,t)$ and $KR1(f,\sigma,2j-1,N,t)$ is denoted as $KR0(2j,t)$ and $KR1(2j,t)$, ($j=1,2,\ldots$). In the plots below, the vertical lines of the grid coincide with the nodes of the grid $\Delta_9$.

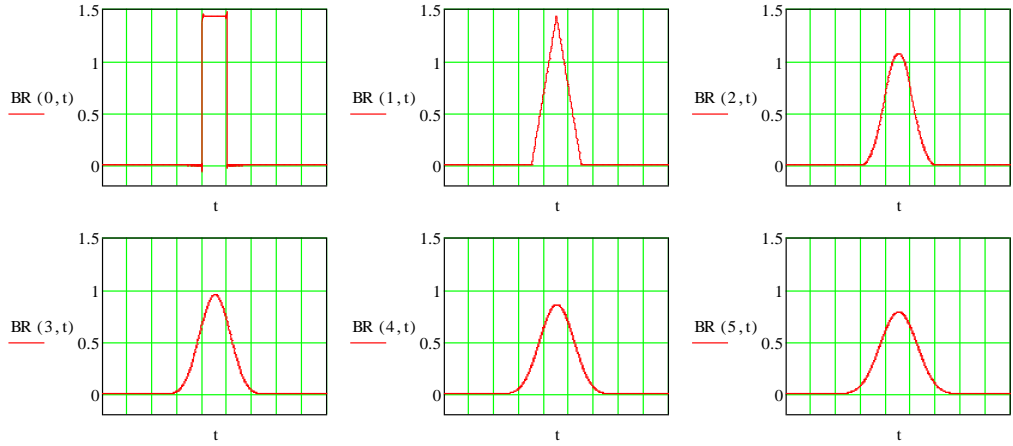

Fig.1. Trigonometric B-splines $BR(r,t)$.

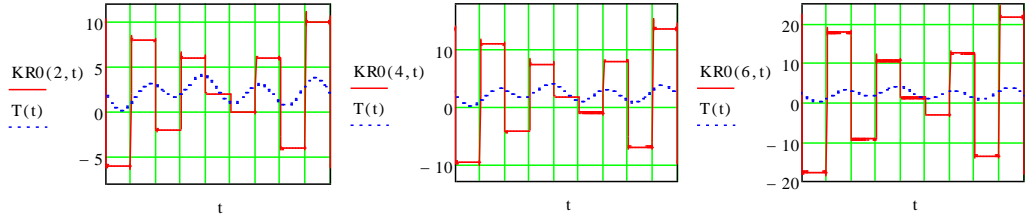

Fig.2. Kernels $KR0(f,\sigma,2j,N,t)$ of interpolation trigonometric splines for even values $j=1,2,3$.

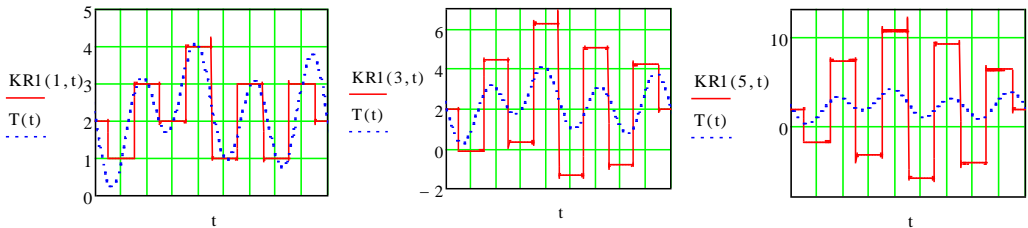

Fig.3. Kernels $KR1(f,\sigma,2j-1,N,t)$ of interpolation trigonometric splines for odd values $j=1,2,3.$.

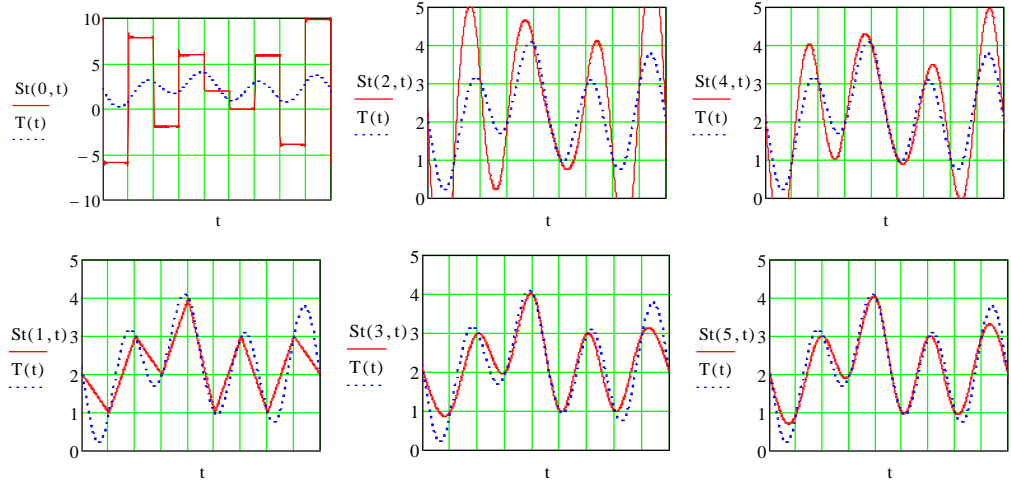

Fig.4. Graphs of interpolation trigonometric splines $St0(f,\sigma,2j,N,t)$ and $St1(f,\sigma,2j+1,N,t)$ for values $j = 0,1,...$ .

There is also another approach to the construction of trigonometric B-splines and kernels of trigonometric interpolation splines. Let us return to formulas (7) and present them in the following way:

$$A_{mn\pm k} = [a_{mN\pm k}\sigma_{mN\pm k}(0,N)]\left[H_k^{-1}(r,N)\sigma_{mN\pm k}(r,N)\right],$$
$$B_{mn\pm k} = [b_{mN\pm k}\sigma_{mN\pm k}(0,N)]\left[H_k^{-1}(r,N)\sigma_{mN\pm k}(r,N)\right] \quad (11)$$
$$(k = 1,2,...,n;\ m = 0,1,...).$$

By applying the same logic as before, we obtain other forms of the trigonometric Riemann spline $BR(r,t)$ of order $r$ ($r = 0,1,...$) and the Riemann kernel $KR0(f,\sigma,r,N,t)$, $KR1(f,\sigma,r,N,t)$ of the interpolation trigonometric spline $St(f,\sigma,r,N,t)$; these forms will be denoted by $BR^*(r,t)$, $KR0^*(f,\sigma,N,t)$ and $KR1^*(f,\sigma,N,t)$ respectively. It is clear that the functions and are represented by the following formulas:

$$BR^*(r,t) = \frac{1}{\pi}\left[\frac{1}{2} + \sum_{k=1}^{\frac{N-1}{2}} \frac{C0(r,k,t)}{Hc(r,k)}\right]; \quad (12)$$

and

$$KR0^*(f,\sigma,N,t) = \frac{a_0}{2} + \sum_{k=1}^{\frac{N-1}{2}}\left[C0(0,k,t)a_k + S0(0,k,t)b_k\right]. \quad (13)$$

$$KR1^*(f,\sigma,N,t) = \frac{a_0}{2} + \sum_{k=1}^{\frac{N-1}{2}}\left[C1(0,k,t)a_k + S1(0,k,t)b_k\right] \quad (13a)$$

Now, (12) demonstrates that the function $BR^*(r,t)$ is obtained from the function $BR(r,t)$ by linear-sumeration methods (see [18]).

As before, we present plots of the functions $BR^*(r,t)$, $KR0^*(f,\sigma,N,t)$ and $KR1^*(f,\sigma,N,t)$ for some values of the parameter $r$.

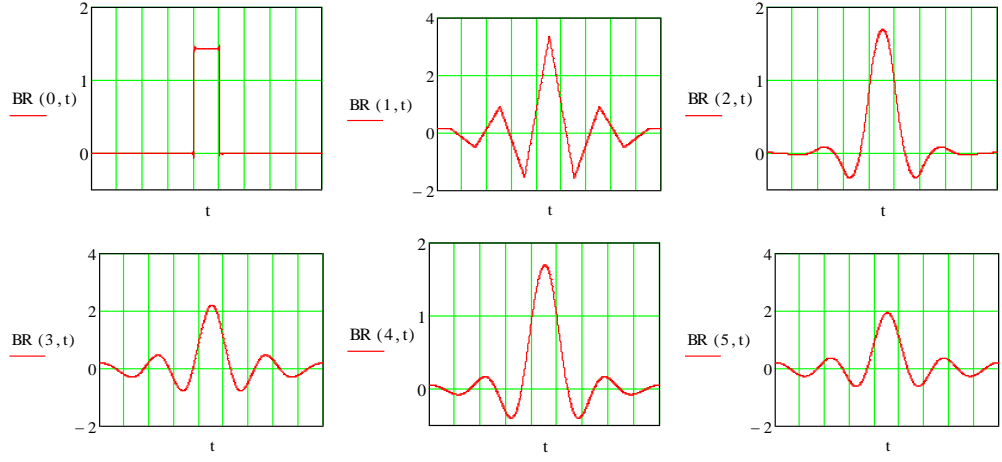

Fig.5. Trigonometric B-splines $BR^*(r,t)$ ($r = 0, 1, \ldots$).

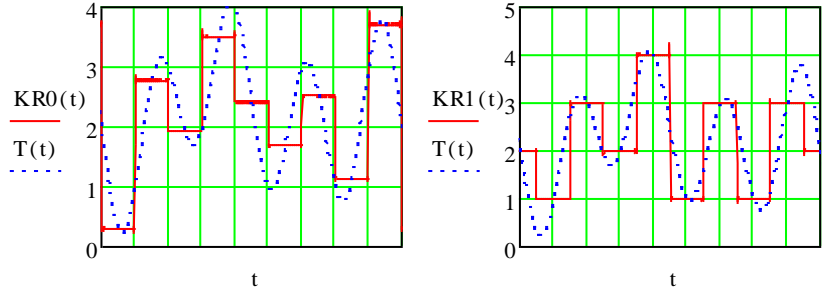

Fig.6. Kernels $KR0^*(f,\sigma,N,t)$ and $KR1^*(f,\sigma,N,t)$ of interpolation trigonometric splines for even and odd values $r$.

Plots of trigonometric interpolation splines $St(f,\sigma,r,N,t)$ are shown in Fig. 4. It is clear that in this case, there also is a correlation

$$St0(f,\sigma,2j,N,t) = \int_0^{2\pi} KR0^*(f,\sigma,N,v) BR^*(2j-1,t-v)dt, \qquad (14)$$

$$St1(f,\sigma,2j-1,N,t) = \int_0^{2\pi} KR1^*(f,\sigma,N,v) BR^*(2(j-1),t-v)dt, \qquad (14a)$$

($j = 1, 2, \ldots$).

Let us have a closer look at the results.

First of all, let us note that the kernels $KR0(f,\sigma,2j,N,t)$, $KR1(f,\sigma,2j-1,N,t)$ and $KR0^*(f,\sigma,N,t)$, $KR1^*(f,\sigma,N,t)$ of trigonometric interpolation splines carry the information about the interpolated function $f(t)$, $t \in [0, 2\pi)$ (in particular, about the grid $\Delta_N$ and the values of the function on this grid). The trigonometric B-splines $BR(r,t)$ and $BR^*(r,t)$ determine the differential properties of the spline $St(f,\sigma,r,N,t)$. Trigonometric B-splines $BR(r,t)$ coincide with polynomial B-splines of the same degree on a segment $[-\pi,\pi]$ with a given grid $\Delta_N$ on it [15]. As a result, the splines $BR(r,t)$ can be called trigonometric B-splines.

It is easy to see that the kernels $KR0(f,\sigma,r,N,t)$ and $KR1(f,\sigma,2j-1,N,t)$ depend on the order $r$, ($r = 1, 2, \ldots$) of the spline; the kernels $KR0^*(f,\sigma,N,t)$ and $KR1^*(f,\sigma,N,t)$ are independent of the order of the spline.

Interpolation trigonometric splines $St(f,\sigma,r,N,t)$ are the result of convolution of trigonometric B-splines $BR(r,t)$ with kernels $KR0(f,\sigma,r,N,t)$ and $KR1(f,\sigma,r,N,t)$ or trigonometric B-splines $BR^*(r,t)$ with kernels $KR0^*(f,\sigma,N,t)$ and $KR1^*(f,\sigma,N,t)$. Since it was shown earlier that interpolation trigonometric splines $St(f,\sigma,2j-1,N,t)$ coincide with periodic polynomial simple interpolation splines of the same degree ($j=1,2,...$), this result can be extended to periodic polynomial simple interpolation splines of odd degrees.

Finally, there is one more matter to pay attention to. It is known [3], [4] that polynomial B-splines are a basis in the space of polynomial splines; this property is widely used (see, e.g., [7], [8]). Trigonometric B-splines $BR(r,t)$ also form a basis in the space of trigonometric splines; this comes from the identity of periodic polynomial and trigonometric interpolation splines. We can assume that trigonometric splines $BR^*(r,t)$ also form a basis in the space of both trigonometric and periodic polynomial splines; however, this fact requires further research.

## Conclusions.

1. The kernels $KR0(f,\sigma,r,N,t)$, $KR1(f,\sigma,r,N,t)$ and $KR0^*(f,\sigma,N,t)$, $KR1^*(f,\sigma,N,t)$ of trigonometric interpolation splines were constructed. Such kernels are the carriers of information about the interpolated function on the grid $\Delta_N$.
2. The classes of trigonometric Riemann B-splines $BR(r,t)$ and $BR^*(r,t)$ were constructed. These B-splines $BR(r,t)$ coincide with polynomial B-splines on the segment $[-\pi,\pi]$.
3. It is shown that interpolation trigonometric splines $St(f,\sigma,r,N,t)$ are the result of convolution of trigonometric B-splines $BR(r,t)$ and $BR^*(r,t)$ with respective kernels $KR0(f,\sigma,r,N,t)$, $KR1(f,\sigma,r,N,t)$ and $KR0^*(f,\sigma,N,t)$, $KR1^*(f,\sigma,N,t)$.
4. In our opinion, the approach proposed in this paper requires further research.